\documentclass[11pt]{article}
\usepackage{amsfonts}
\usepackage{color}
\usepackage{amsmath,amssymb,comment}
\usepackage{verbatim}

\newtheorem{theo}{Theorem}[section]
\newtheorem{lem}[theo]{Lemma}

\newtheorem{defi}[theo]{Definition}

\newcommand{\mysection}[1]{\section{#1} \setcounter{equation}{0}}
    \makeatletter
\def\@fnsymbol#1{\ensuremath{\ifcase#1\or *\or \ddagger\or
   \mathsection\or \mathparagraph\or \|\or **\or \dagger\dagger
   \or \ddagger\ddagger \else\@ctrerr\fi}}
    \makeatother
\newcommand{\proof}{{\sc Proof.} \quad}
\newcommand{\proofc}{{\sc Proof} \ }
\newcommand{\be}{\begin{equation} \label}
\newcommand{\ee}{\end{equation}}
\newcommand{\bea}{\begin{eqnarray}\label}
\newcommand{\eea}{\end{eqnarray}}
\newcommand{\bas}{\begin{eqnarray*}}
\newcommand{\eas}{\end{eqnarray*}}
\newcommand{\bit}{\begin{itemize}}
\newcommand{\eit}{\end{itemize}}
\newcommand{\qed}{\hfill$\Box$ \vskip.2cm}
\newcommand{\nn}{\nonumber}
\newcommand{\R}{\mathbb{R}}

\newcommand{\pO}{\partial\Omega}

\newcommand{\io}{\int_\Omega}

\newcommand{\na}{\nabla}
\newcommand{\Del}{\Delta}
\newcommand{\del}{\delta}
\newcommand{\al}{\alpha}

\newcommand{\sig}{\sigma}

\newcommand{\pa}{\partial}
\newcommand{\bom}{\overline{\Omega}}
\newcommand{\Om}{\Omega}

\newcommand{\mint}{- \hspace*{-3.3mm} \int}

\newcommand{\ov}{\overline}
\newcommand{\un}{\underline}
\newcommand{\wh}{\widehat}

\newcommand{\hs}{\hspace*}
\newcommand{\sm}{\setminus}

\newcommand{\vp}{\varphi}
\newcommand{\lbal}{\left\{ \begin{array}{l}}
\newcommand{\lball}{\left\{ \begin{array}{ll}}
\newcommand{\ear}{\end{array} \right.}
\newcommand{\eeeear}{\end{array}}

\newcommand{\abs}{\\[5pt]}

\newcommand{\tm}{T_{max}}

\newcommand{\uU}{\un{U}}
\newcommand{\oU}{\ov{U}}

\newcommand{\uW}{\un{W}}
\newcommand{\oW}{\ov{W}}

\newcommand{\hU}{\wh{U}}
\newcommand{\hW}{\wh{W}}

\newcommand{\muu}{\mu^{(u)}}
\newcommand{\muw}{\mu^{(w)}}
\newcommand{\mus}{\mu_\star}
\newcommand{\muS}{\mu^\star}
\newcommand{\Ms}{M_\star}
\newcommand{\MS}{M^\star}

\newcommand{\parab}{{\mathcal{P}}}
\newcommand{\qarab}{{\mathcal{Q}}}

\newcommand{\gams}{\gamma_\star}

\newcommand{\dels}{\del_\star}
\newcommand{\delss}{\del_{\star\star}}

\newcommand{\ys}{y_\star}

\newcommand{\sst}{s_\star}
\newcommand{\ssst}{s_{\star\star}}

\newcommand{\ths}{\theta_\star}
\newcommand{\thss}{\theta_{\star\star}}
\newcommand{\Ts}{T_\star}
\newcommand{\TS}{T^\star}
\begin{document}
\enlargethispage{10mm}
\title{A critical nonlinearity for blow-up in a higher-dimensional chemotaxis system with indirect signal production}
\author
{
Yiheng Zhao\footnote{zhaoyh2002@sjtu.edu.cn}\\
{\small School of Mathematical Sciences, Shanghai Jiao Tong University,}\\
{\small Shanghai 200240, P.R.~China} }
\date{}
\maketitle
\begin{abstract}
\noindent
  The Neumann problem in balls $\Om\subset\R^n$, $n\in\{3,4\}$, for the chemotaxis system
  \bas
	\left\{ \begin{array}{ll}	
	u_t = \Del u - \na \cdot (u\na v), \\[1mm]
	0 = \Del v - \muw(t) + w, \qquad \muw(t) = \mint_\Om w,\\[1mm]
	w_t = \Del w - w + f(u),
	\end{array} \right.
  \eas
  is considered. Under the assumption that $f\in C^1([0,\infty))$ is such that $f(\xi) \ge k\xi^\sig$
  for all $\xi\ge 0$ and some $k>0$ and $\sig>\frac{4}{n}$, it is shown that finite-time blow-up
  occurs for some radially symmetric solutions.\abs
\end{abstract}
\newpage
\section{Introduction}\label{intro}
Chemotaxis systems including indirect signal production mechanisms on the one hand appear to provide more realistic
descriptions of cross-diffusive migration than simple Keller-Segel systems in various application contexts
(\cite{shi_jmb2021}, \cite{shi_zamp2021}, \cite{macfarlane}, \cite{MPB});
in line with this, accordingly extended models have been in the focus of an increasing number of analytical investigations.
Besides cases in which the respective attractant production is mediated by a non-diffusible quantity (\cite{fuest_m2},
\cite{laurencot_stinner_sima}, \cite{hu_tao_M3AS}, \cite{taowin_JEMS}), a considerable literature has been concerned with related situations
in which the corresponding additional component responsible for signal production diffuses in space
(\cite{fujie_senba}, \cite{ding_wang}, \cite{wei_wang}, \cite{wu}).\abs
Indeed, several studies on particular representatives of the general model class
\be{00}
	\left\{ \begin{array}{ll}	
	u_t = \na \cdot (D(u)\na u) - \na \cdot (S(u)\na v), \\[1mm]
	v_t = \Del v - v + w, \\[1mm]
	w_t = \Del w - w + f(u),
	\end{array} \right.
\ee
have revealed a significant strengthening of the tendency toward boundedness in comparison to two-component
relatives accounting for direct signal production.
In the most prototypical case when here $D\equiv 1, S\equiv id$ and $f\equiv id$, for instance,
blow-up phenomena can be ruled out not only in two- but even in three-dimensional boundary value problems (\cite{fujie_senba});
also in the presence of more general and suitably regular migration rates $D$ and $S$, results on global solvability and
boundedness addressing (\ref{00}) require significantly weaker assupmtions on the large-density
asymptotics of $S$, relative to that of $D$, than those known to be critical for corresponding
two-component quasilinear systems (\cite{ding_wang}, \cite{cieslak_stinner_JDE2012}, \cite{win_collapse}).
For some further results concerned with variants of (\ref{00}), we refer to \cite{zheng_m3as2022}
and \cite{li_xiang_jde2023}.\abs
Accordingly, detections of singularity formation in models of this and related types can yet only rarely
be found in the analytical literature.
When again $S\equiv f\equiv id$, for instance, in line with the above a necessary condition for blow-up to be possible
at all is that the diffusion rate $D(u)$ undergoes a suitably substantial decay at large population densities $u$,
and a corresponding quantification has recently been achieved for a parabolic-elliptic simplification of (\ref{00})
with $f\equiv id$, actually throughout large classes of fairly general power-type ingredients $D$ and $S$
(\cite{taowin_291} and \cite{taowin_292}).
Beyond this, rigorous findings on the occurrence of unboundedness phenomena in chemotaxis models with indirect
signal production seem to have concentrated on 		%less strongly dissipation-influenced
scenarios involving non-diffusible attractant producers (\cite{laurencot_stinner_sima}, \cite{taowin_JEMS}, \cite{fuest_m2}).
This may be viewed as confirming a substantially stronger dissipative character of (\ref{00}) in comparison to
the latter model class; correspondingly, the potential for a derivation of blow-up results on the basis
of energy inequalities involving Lyapunov functionals unbounded from below, as known from the analysis of classical
Keller-Segel systems (\cite{cieslak_stinner_JDE2012}, \cite{win_JMPA}), seems limited due to the presence of lower bounds
for energy functionals even in physically meaningful subcases of (\ref{00}) for which they have been found (\cite{fujie_senba}).\abs
{\bf Main results.} \quad
The intention of this manuscript is to address this methodological problem, and to develop a technique which is capable
of extending the yet small selection of 		%cases in which singularity formation can be detected in
chemotaxis systems involving diffusion-influenced indirect signal production mechanisms
in which singularity formation can be detected.
This will be pursued in the framework of the problem
\be{0}
	\left\{ \begin{array}{ll}	
	u_t = \Del u - \na \cdot (u\na v),
	\qquad & x\in\Om, \ t>0, \\[1mm]
	0 = \Del v - \muw(t) + w,
	\qquad \muw(t) = \mint_\Om w,
	\qquad & x\in\Om, \ t>0, \\[1mm]
	w_t = \Del w - w + f(u)
	\qquad & x\in\Om, \ t>0, \\[1mm]
	\frac{\pa u}{\pa\nu}=\frac{\pa v}{\pa\nu}=\frac{\pa w}{\pa\nu}=0,
	\qquad & x\in\pO, \ t>0, \\[1mm]
	u(x,0)=u_0(x), \quad w(x,0)=w_0(x),
	\qquad & x\in\Om,
	\end{array} \right.
\ee
which can be viewed as a straightforward parabolic-elliptic-parabolic simplification of (\ref{00}), justifiable in cases
when the signal concentration $v$ diffuses suitably fast (\cite{JL}), in the special situation when
$D\equiv 1$ and $S\equiv id$, but when now the signal production rate $f$ may depend nonlinearly on the population density $u$.
For the corresponding version of (\ref{00}) it is known that in any smoothly bounded domain $\Om\subset\R^n$ with $n\ge 2$,
the assumption that here
\be{f1}
0\le f\in C^1([0,\infty))
\ee
and
\be{f01}
	f(\xi) \le K \xi^\sig
	\qquad \mbox{for all } \xi>1,
\ee
with some $K>0$ and some
\be{f02}
	\sig<\frac{4}{n}
	\quad \mbox{when } n\ge 2;
    \qquad\mbox{or }
    \quad
    \sig<3
    \quad \mbox{when } n=1
\ee
is sufficient to ensure that all suitably regular nonnegative initial data evolve into global bounded solutions
(\cite{cao_aml2024}, \cite{tao_zhang_aml2024}); in fact, by simple adaptation of the arguments therein this statement can readily be carried over
to (\ref{0}).\abs
In order to supplement this by a result on the occurrence of blow-up under essentially optimal assumptions on the growth
of $f$, we shall further develop an approach recently presented in \cite{taowin_291}, based on the observation
that in settings of radially symmetric solutions in balls, the corresponding cumulated densities associated with $u$ and $w$
satisfy a parabolic system which is cooperative and hence admits a comparison principle (Lemma 3.1 and Lemma 3.2).
The design of suitable exploding subsolutions thereof, inspired by a corresponding construction in
\cite{taowin_291}, indeed will enable us to reveal the presence of blow-up phenomena
in three- and four-dimensional versions of (\ref{0})
within ranges of ingredients $f$ which, in view of the above implication on (\ref{f01}) and (\ref{f02}) on boundedness,
apparently cannot be substantially extended.
More precisely, our main result is the following.
\begin{theo}\label{theo13}
  Let $n\in \{3, 4\}, R>0$ and $\Om=B_R\equiv B_R(0)\subset\R^n$, and suppose that $f$
  satisfies (\ref{f1}) as well as
  \be{f2}
  f(u)\ge k u^\sig
  \qquad\mbox{for all } u\ge 0
  \ee
 with some $k>0$ and
    \be{sig}
	\sig>\frac{4}{n}.
  \ee
  Then one can find some initial data $(u_0, v_0)$ fulfilling
  \be{ir}
	u_0\in W^{1,\infty}(\Om)
	\mbox{ and }
	w_0\in W^{1,\infty}(\Om)
	\mbox{ are radially symmetric and nonnegative}
  \ee
  and there exists $\TS>0$ such that
  the corresponding classical solution $(u,v,w)$ of (\ref{0}) blows up before time $\TS$ in the following sense:
  There exists $\tm\in (0,\TS)$ such that
  \be{13.3}
	\limsup_{t\nearrow\tm} \|u(\cdot,t)\|_{L^\infty(\Om)}=\infty.
  \ee
\end{theo}
From Theorem \ref{theo13} together with the boundedness results in \cite{tao_zhang_aml2024} and \cite{cao_aml2024}
we infer that $\sig=\frac{4}{n}$ ($n=3, 4$) is the critical exponent for blow-up to (\ref{0}).
In comparison with the corresponding two-component Keller-Segel systems with nonlinear signal production, we
reminisce that $\sigma=\frac{2}{n}$ ($n\ge 1$) is the critical blow-up exponent for the latter systems:
Namely, if $f$ satisfies (\ref{f1}), (\ref{f2}) and $\sig>\frac{2}{n}$, then the solutions blow up in finite time for suitably large
initial-data (\cite{winkler_non2018}); whereas the solutions
exist globally and are bounded for any given smooth initial data whenever $f$ fulfills (\ref{f1}),
(\ref{f01}) and $\sig<\frac{2}{n}$ (\cite{liu_tao2016}).
\mysection{Local existence}
We start from recalling a statement concerning the local existence and an extensibility criterion.
\begin{lem}\label{lem_loc}
  Let $n\ge 1$ and $\Om\subset\R^n$ be a bounded domain with smooth boundary,
  and assume (\ref{f1}) and (\ref{ir}).
  Then there exist $\tm\in (0,\infty]$ and a uniquely determined triple of functions
  \be{reg_loc}
	\lbal
	u\in C^0(\bom\times [0,\tm)) \cap C^{2,1}(\bom\times (0,\tm)), \\[1mm]
	v \in C^{2,0}(\bom\times (0,\tm))
	\qquad \mbox{and} \\[1mm]
	w\in C^0(\bom\times [0,\tm)) \cap C^{2,1}(\bom\times (0,\tm))
	\ear
  \ee
  with the properties that $u\ge 0$ and $w\ge 0$ in $\Om\times (0,\tm)$, that $\io v(\cdot,t)=0$ for all $t\in (0,\tm)$,
  that (\ref{0}) is satisfied in the classical sense in $\Om\times (0,\tm)$, and that
  \be{ext}
	\mbox{if $\tm<\infty$, \qquad then \qquad}
	\limsup_{t\nearrow\tm} \|u(\cdot,t)\|_{L^\infty(\Om)} =\infty.
  \ee
  Moreover,
  \be{mass}
	\io u(\cdot,t) dx = \io u_0 dx
	\qquad \mbox{for all } t\in (0,\tm).
  \ee
\end{lem}
\proof
  The local existence as well as the extensibility criterion (\ref{ext}) can be achieved by adapting
  fixed-point approaches well-established in the context of chemotaxis systems
  (cf., e.g., \cite{ding_wang} for some close relative). The mass identity (\ref{mass}) simply results
  from an integration of the first equation in (\ref{0}).
\qed
\mysection{Finite-time blow-up} \label{sect_BU}
\subsection{A cooperative parabolic system for mass functions} \label{sect_coop}
Motivated by an idea that the construction of radial blow-up solutions in Kelle-Segel type systems
can be achieved via analysis of associated cumulated densities (\cite{JL}, \cite{biler1998} and \cite{taowin_291}),
we first transform the original system (\ref{0}) into a parabolic system
satisfied by cumulated densities in radially symmetric setting. Here we point out that in order to successfully deal with the nonlinear
zero-order production term $f(u)$, we will need the spatial-dimension restriction that $n\le 4$.
\begin{lem}\label{lem11}
  Let $n\le 4$ and $R>0$, assume (\ref{f1}) and (\ref{f2}) as well as (\ref{sig}), and let $\Ms>0$ and $\MS>\Ms$.
  Then, whenever (\ref{ir}) holds  with
  \be{11.1}
	\Ms \le \io u_0 dx \le \MS,
	\qquad
	\Ms \le \io w_0 dx \le \MS
    \qquad \mbox{and} \qquad
    \|w_0\|_{L^\infty(\Om)} \le \frac{\MS}{|\Om|},
  \ee
  there exists $T_0=T_0^{(u_0, w_0)}>0$ such that
   the corresponding solution $(u,v,w)$ of (\ref{0}) has the property that setting
  \be{U}
	U(s,t):=\int_0^{s^\frac{1}{n}} \rho^{n-1} u(\rho,t) d\rho,
	\qquad s\in [0,R^n], \ t\in [0,\tm),
  \ee
  and
  \be{W}
	W(s,t):=\int_0^{s^\frac{1}{n}} \rho^{n-1} w(\rho,t) d\rho,
	\qquad s\in [0,R^n], \ t\in [0,\tm),
  \ee
  we obtain functions $U$ and $W$ which belong to $C^1([0,R^n]\times [0,\tm)) \cap C^{2,1}((0,R^n) \times (0,\tm))$
  and satisfy $U_s\ge 0$ and $W_s\ge 0$ as well as
  \be{11.2}
	U(0,t)=W(0,t)=0,
	\quad
	U(R^n,t) \ge \frac{\mus R^n}{n}
	\quad \mbox{and} \quad
	W(R^n,t) \ge \frac{\mus R^n}{n}
	\qquad \mbox{for all $t\in (0,\tm)\cap (0,T_0)$,}
  \ee
  and for which we have
  \be{11.3}
	\parab^{(\muS)} [U,W](s,t) \ge 0
	\qquad \mbox{for all $s\in (0,R^n)$ and $t\in (0,\tm)\cap (0,T_0)$}
  \ee
  as well as
  \be{11.4}
	\qarab [U,W](s,t) \ge 0
	\qquad \mbox{for all $s\in (0,R^n)$ and $t\in (0,\tm)\cap (0,T_0)$,}
  \ee
  where
  \be{11.5}
	\mus:=\frac{\Ms}{2|\Om|}
	\qquad \mbox{and} \qquad
	\muS:=\frac{2\MS}{|\Om|}.
  \ee
  Here and below, for $T>0$ and functions $\vp$ and $\psi$ from $C^1([0,R^n]\times [0,T))$ which satisfy $\vp_s \ge 0$ on
  $(0,R^n)\times (0,T)$ and are such that
  $\vp(\cdot,t) \in W^{2,\infty}_{loc}((0,R^n))$ and
  $\psi(\cdot,t) \in W^{2,\infty}_{loc}((0,R^n))$ for all $t\in (0,T)$, we let
  \be{P}
	\parab^{(\muS)}[\vp,\psi](s,t)
	:= \vp_t - n^2 s^{2-\frac{2}{n}} \vp_{ss}
	- n\vp_s \cdot \Big(\psi-\frac{\muS s}{n}\Big)
   \ee
  and
  \be{Q}
	\qarab[\vp,\psi](s,t)
	:= \psi_t - n^2 s^{2-\frac{2}{n}} \psi_{ss} + \psi -K s^{1-\sig} \vp^\sig
  \ee
  for $t\in (0,T)$ and a.e.~$s\in (0,R^n)$, where $K :=kn^{\sig-1}$.
\end{lem}
\proof The mass identity (\ref{mass}) along with the first inequality in (\ref{11.1}) entails that
 \be{11.6}
 \muu(t)\equiv \mint_\Om u(\cdot,t) dx =\mint_\Om u_0 dx
 \ge \frac{\Ms}{|\Om|} \ge \frac{\Ms}{2|\Om|}
 \qquad \mbox{for all } t\in (0,\tm),
 \ee whence an integration of the third equation in (\ref{0}) together with the nonnegativity of $f(u)$
 yields
 \bas
 \frac{d}{dt} \muw(t) + \muw(t)\ge 0
 \qquad \mbox{for all } t\in (0,\tm)
 \eas that implies
 \be{11.8}
 \muw(t)\ge e^{-t}\muw(0) \ge \frac{1}{2} \muw(0) \ge \frac{\Ms}{2|\Om|}
 \qquad \mbox{for all } t\in (0,\tm)\cap (0,\ln 2)
 \ee due to the second inequality in (\ref{11.1}), where $\muw(t)\equiv \mint_\Om w(\cdot,t) dx$.
 On the other hand, it follows from the continuity of the solution $(u, v, w)$ that there exists $T_0=T_0^{(u_0, w_0)} \in (0, \ln 2]$ such that
  \bas
   \|w(\cdot, t)\|_{L^\infty(\Om)} \le 2 \|w_0\|_{L^\infty(\Om)}
   \qquad \mbox{for all } t\in (0,\tm)\cap (0,T_0)
  \eas  and thus thanks to the third inequality in (\ref{11.1}),
  \be{11.9}
   \muw(t)=\frac{1}{|\Om|} \io w
   \le \|w(\cdot, t)\|_{L^\infty(\Om)}
   \le 2 \|w_0\|_{L^\infty(\Om)}
   \le \frac{2\MS}{|\Om|}
   \qquad \mbox{for all } t\in (0,\tm)\cap (0,T_0).
  \ee
  Since a direct computation based on (\ref{0}) shows that (\ref{U}) and (\ref{W}) indeed define functions
  $U$ and $W$ which possess the claimed regularity and monotonicity properties, and which moreover satisfy
  \be{11.10}
	U(0,t)=W(0,t)=0,
	\quad
	U(R^n,t) = \frac{\muu(t) R^n}{n}
	\quad \mbox{and} \quad
	W(R^n,t) = \frac{\muw(t) R^n}{n}
	\qquad \mbox{for all } t\in (0,\tm)
  \ee
  as well as
   \be{11.11}
	U_t= n^2 s^{2-\frac{2}{n}} U_{ss} + nU_s \cdot \Big(W-\frac{\muw(t) s}{n}\Big)
	\qquad \mbox{in } (0,R^n)\times (0,\tm)
  \ee
  and
  \be{11.12}
	W_t= n^2 s^{2-\frac{2}{n}} W_{ss} -W +\frac{1}{n} \int_0^s f(n U_\xi(\xi, t)) d\xi
	\qquad \mbox{in } (0,R^n)\times (0,\tm).
  \ee
  From (\ref{11.10}), (\ref{11.6}) and (\ref{11.8}) we obtain (\ref{11.2}), while (\ref{11.3})
  is a consequence of (\ref{11.11}) and (\ref{11.9}) along with the nonnegativity of $U_s$.
  In order to verify (\ref{11.4}), we first note that
 $\sig>1$ due to (\ref{sig}) together with our assumption that $n\le 4$, and thus we can use the H\"{o}lder inequality
 to obtain
 \bas
 U^\sig(s, t) =\bigg(\int_0^s U_\xi(\xi, t) d\xi\bigg)^\sig \le s^{\sig-1} \int_0^s U_\xi^\sig(\xi, t) d\xi
 \eas thanks to the nonnegativity of $U_\xi$,
 which along with (\ref{f2}) leads to
 \bas
 \frac{1}{n} \int_0^s f(n U_\xi(\xi, t)) d\xi
 \ge k n^{\sig-1} \int_0^s U_\xi^\sig(\xi, t) d\xi
 \ge K s^{1-\sig}U^\sig(s, t)
  \eas
  with $K:=k n^{\sig-1}$.
 Combining (\ref{11.12}) with this we finally arrive at (\ref{11.4}).
\qed
A cooperative feature of the parabolic system associated to (\ref{P}) and (\ref{Q}) guarantees
that a comparison principle holds, which is a corner stone of the construction of a blow-up sub-solution.
\begin{lem}\label{lem_CP}
  Let $n\le 4$ and $R>0$, assume (\ref{sig}),
  and suppose that for some $T>0$, the functions
  $\uU,\oU,\uW$ and $\oW$ belong to $C^1([0,R^n] \times [0,T))$ and have the properties that
  \bas
  \{\uU(\cdot,t),\oU(\cdot,t),\uW(\cdot,t),\oW(\cdot,t)\} \subset W^{2,\infty}_{loc}((0,R^n))
  \qquad \mbox{for all } t\in (0,T),
  \eas
  and that $\uU_s$ and $\oU_s$ are nonnegative on $(0,R^n)\times (0,T)$.
  If with $\parab^{(\muS)}$ and $\qarab=\qarab^{(\muS)}$ as in (\ref{P}) and (\ref{Q}) we have
  \be{C1}
	\parab^{(\muS)}[\uU,\uW](s,t) \le 0
	\quad \mbox{and} \quad
	\parab^{(\muS)}[\oU,\oW](s,t) \ge 0
	\qquad \mbox{for all $t\in (0,T)$ and a.e.~} s\in (0,R^n)
  \ee
  as well as
  \be{C2}
	\qarab[\uU,\uW](s,t) \le 0
	\quad \mbox{and} \quad
	\qarab[\oU,\oW](s,t) \ge 0
	\qquad \mbox{for all $t\in (0,T)$ and a.e.~} s\in (0,R^n),
  \ee
  and if furthermore
  \be{C3}
	\uU(0,t) \le \oU(0,t),
	\quad
	\uU(R^n,t) \le \oU(R^n,t),
	\quad
	\uW(0,t) \le \oW(0,t)
	\quad \mbox{and} \quad
	\uW(R^n,t) \le \oW(R^n,t)
%	\qquad \mbox{for all } t\in [0,T)
  \ee
  for all $t\in [0,T)$
  as well as
  \be{C4}
	\uU(s,0) \le \oU(s,0)
	\quad \mbox{and} \quad
	\uW(s,0)\le \oW(s,0)
	\qquad \mbox{for all } s\in [0,R^n],
  \ee
  then
  \be{C5}
	\uU(s,t)\le \oU(s,t)
	\qquad \mbox{for all $s\in [0,R^n]$ and } t\in [0,T).
  \ee
\end{lem}
\proof Since for any given $t_0\in (0, T)$, in view of the fact that $\sig> 1$ due to (\ref{sig}) and $n\le 4$ and according to our regularity assumption on $\uU$ and $\oU$,
the Lagrange intermediate value theorem guarantees the existence of some function $A(s, t_0)$ and a constant $c(t_0)>0$
such that
\bas
\uU^\sig(s, t_0) -\oU^\sig(s, t_0) =A(s, t_0)(\uU (s, t_0) -\oU (s, t_0))
\qquad\mbox{for all } s\in [0, R^n]
\eas with $|A(s, t_0)|\le c(t_0)$, the proof is then similar to \cite[Lemma 3.2]{taowin_291} and therefore we refrain
 us from repeating the details here.
\qed
\subsection{Construction of blow-up subsolutions}
The blow-up subsolutions will involve two crucial parameters $\alpha$ and $\beta$. However, under the
assumption that $\sig>\frac{4}{n}$, the choice of these two parameters satisfying (\ref{30.2})
below will require another spatial-dimension restriction that $n\ge 3$.
\begin{lem}\label{lem30}
Let $n\ge 3$, and assume
\be{30.1}
\sig>\frac{4}{n}.
\ee
Then there exist $\al \in (0, 1)$ and $\beta\in (0, 1)$ such that
\be{30.2}
1+\frac{2}{n} -\sig +\al \sig <\beta < 1-\frac{2}{n}.
\ee
\end{lem}
\proof  Since (\ref{30.1}) ensures that
\bas
1+\frac{2}{n} -\sig < 1-\frac{2}{n},
\eas we can take $\al =\al(\sig)\in (0, 1)$ sufficiently small such that
\bas
1+\frac{2}{n} -\sig +\al \sig  < 1-\frac{2}{n}
\eas still holds. Due to $n\ge 3$, we have $1-\frac{2}{n}\ge \frac{1}{3}>0$, and thus
we can fix some $\beta\in (0, 1)$ fulfilling (\ref{30.2}).
\qed
Although the essential part of the blow-up subsolutions will take the same form as in \cite{taowin_291},
the choices of parameters $\alpha, \beta$ and parameter function $y(t)$ here are different from those in \cite{taowin_291}.
\begin{lem}\label{lem1}
  Let $n\in \{3, 4\}$, $R>0$ and $\mus>0$, let $\al \in (0, 1)$ be provided by Lemma \ref{lem30},
  and let
  \be{a}
	a\equiv a^{(\mus)} := \frac{\mus R^n}{n e^\frac{1}{e} (R^n+1)}.
  \ee
   Then for all $T>0$ and any $y\in C^1([0,T))$ with $y(t)>\frac{1}{R^n}$ for all $t\in (0,T)$, setting
  \be{hU}
	\hU(s,t) \equiv \hU^{(\mus,\al,y)} :=
	\lball
	a y^{1-\al}(t) s,
	\qquad & t\in [0,T), \ s\in \Big[ 0, \frac{1}{y(t)}\Big], \\[1mm]
	\al^{-\al} a \cdot \Big(s-\frac{1-\al}{y(t)} \Big)^\al,
	\qquad & t\in [0,T), \ s\in \Big( \frac{1}{y(t)} , R^n \Big],
	\ear
  \ee
  defines a function $\hU\in C^1([0,R^n] \times [0,T)) \cap C^0([0,T);W^{2,\infty}((0,R^n))$ which satisfies
  \be{1.1}
	\hU(0,t)=0
	\quad \mbox{and} \quad
	\hU(R^n,t) \le \frac{\mus R^n}{n}
	\qquad \mbox{for all } t\in (0,T),
  \ee
  and for which we have $\hU(\cdot,t) \in C^2([0,R^n]\sm \{\frac{1}{y(t)}\}$ for all $t\in (0,T)$, with
  \be{1.2}
	\hU_t(s,t) = \lball
	(1-\al) a y^{-\al}(t) y'(t) s,
	\qquad t\in (0,T), \ s\in (0,\frac{1}{y(t)}), \\[1mm]
	\al^{1-\al} (1-\al) a \cdot \Big(s-\frac{1-\al}{y(t)} \Big)^{\al-1} \cdot \frac{y'(t)}{y^2(t)},
	\qquad & t\in (0,T), \ s\in \Big( \frac{1}{y(t)} , R^n \Big),
	\ear
  \ee
  and
  \be{1.3}
	\hU_s(s,t) = \lball
	a y^{1-\al}(t),
	\qquad t\in (0,T), \ s\in (0,\frac{1}{y(t)}), \\[1mm]
	\al^{1-\al} a \cdot \Big(s-\frac{1-\al}{y(t)} \Big)^{\al-1},
	\qquad & t\in (0,T), \ s\in \Big( \frac{1}{y(t)} , R^n \Big),
	\ear
  \ee
  as well as
  \be{1.4}
	\hU_{ss}(s,t) = \lball
	0,
	\qquad t\in (0,T), \ s\in (0,\frac{1}{y(t)}), \\[1mm]
	- \al^{1-\al} (1-\al) a \cdot \Big(s-\frac{1-\al}{y(t)} \Big)^{\al-2},
	\qquad & t\in (0,T), \ s\in \Big( \frac{1}{y(t)} , R^n \Big).
	\ear
  \ee
\end{lem}
\proof
  The claimed regularity properties and the identities in (\ref{1.2})-(\ref{1.4})
  as well as the inequality in (\ref{1.1}) can be directly verified on the basis of (\ref{hU}),
  (\ref{a}) and the restriction $\al\in (0, 1)$.
\qed
Similarly, we also have:
\begin{lem}\label{lem2}
  Let $n\in \{3, 4\}$, $R>0$, $\mus>0$, and let $\beta\in (0,1)$ be taken from Lemma \ref{lem30}.
  Then whenever $T>0$ and $y\in C^1([0,T))$ is such that $y(t)>\frac{1}{R^n}$ for all $t\in (0,T)$, taking $a=a^{(\mus)}$
  from (\ref{a}) and writing
  \be{hW}
	\hW(s,t) \equiv \hW^{(\mus,\beta,y)}(s,t):= \lball
	a y^{1-\beta}(t) s,
	\qquad & t\in [0,T), \ s\in \Big[ 0, \frac{1}{y(t)}\Big], \\[1mm]
	\beta^{-\beta} a \cdot \Big(s-\frac{1-\beta}{y(t)} \Big)^\beta,
	\qquad & t\in [0,T), \ s\in \Big( \frac{1}{y(t)} , R^n \Big],
	\ear
  \ee
  we obtain an element $\hW$ of $C^1([0,R^n] \times [0,T)) \cap C^0([0,T);W^{2,\infty}((0,R^n))$ with the properties that
  \be{2.1}
	\hW(0,t)=0
	\quad \mbox{and} \quad
	\hW(R^n,t) \le \frac{\mus R^n}{n}
	\qquad \mbox{for all } t\in (0,T),
  \ee
  that
  \be{2.2}
	\hW_t(s,t) = \lball
	(1-\beta) a y^{-\beta}(t) y'(t) s,
	\qquad t\in (0,T), \ s\in (0,\frac{1}{y(t)}), \\[1mm]
	\beta^{1-\beta} (1-\beta) a \cdot \Big(s-\frac{1-\beta}{y(t)} \Big)^{\beta-1} \cdot \frac{y'(t)}{y^2(t)},
	\qquad & t\in (0,T), \ s\in \Big( \frac{1}{y(t)} , R^n \Big),
	\ear
  \ee
  that
  \be{2.3}
	\hW_s(s,t) = \lball
	a y^{1-\beta}(t),
	\qquad t\in (0,T), \ s\in (0,\frac{1}{y(t)}), \\[1mm]
	\beta^{1-\beta} a \cdot \Big(s-\frac{1-\beta}{y(t)} \Big)^{\beta-1},
	\qquad & t\in (0,T), \ s\in \Big( \frac{1}{y(t)} , R^n \Big),
	\ear
  \ee
  and that $\hW(\cdot,t) \in C^2([0,R^n]\sm \{\frac{1}{y(t)}\}$ for all $t\in (0,T)$, with
  \be{2.4}
	\hW_{ss}(s,t) = \lball
	0,
	\qquad t\in (0,T), \ s\in (0,\frac{1}{y(t)}), \\[1mm]
	- \beta^{1-\beta} (1-\beta) a \cdot \Big(s-\frac{1-\beta}{y(t)} \Big)^{\beta-2},
	\qquad & t\in (0,T), \ s\in \Big( \frac{1}{y(t)} , R^n \Big).
	\ear
  \ee
\end{lem}
\proof
  The proof is the same as that of Lemma \ref{lem1}.
\qed
Now we can define the following family of functions that act as candidates for subsolutions.
\begin{defi}\label{duw}
  Let $n\in \{3, 4\}$, $R>0$, $\mus>0$ and $\theta>0$,
  let $\al\in (0,1)$ and $\beta\in (0,1)$ be chosen by Lemma \ref{lem30},
  and let $a=a^{(\mus)}$ be as in (\ref{a}). Given $T>0$ and
  $y\in C^1([0,T))$ such that $y(t)>\frac{1}{R^n}$ for all $t\in [0,T)$, we then let
  \be{uU}
	\uU(s,t) \equiv \uU^{(\mus,\al,\theta,y)}(s,t):=e^{-\theta t} \hU(s,t),
	\qquad s\in [0,R^n], \ t\in [0,T),
  \ee
  and
  \be{uW}
	\uW(s,t) \equiv \uW^{(\mus,\beta,\theta,y)}(s,t):=e^{-\theta t} \hW(s,t),
	\qquad s\in [0,R^n], \ t\in [0,T),
  \ee
  where $\hU=\hU^{(\mus,\al,y)}$ and
  $\hW=\hW^{(\mus,\beta,y)}$ are taken from (\ref{hU}) and (\ref{hW}).
\end{defi}
From now on,  we take $\sig>\frac{4}{n}$ with $n\in \{3, 4\}$ and then fix $\al \in (0, 1)$ and $\beta\in (0, 1)$
satisfying (\ref{30.2}).
\subsection{Subsolution properties: inner region}
We begin with verifying that $(\uU,\uW)$ satisfies $\parab^{(\muS)} [\uU,\uW](s,t) \le 0$
near $s=0$.
\begin{lem}\label{lem5}
  Let $n\in \{3, 4\}$ and $R>0$, and assume (\ref{sig}).
  %, and let $\al\in (0, 1)$ and $\beta \in (0, 1)$ from Lemma \ref{lem30}.
  Then there exists $\dels\in (0, 1)$ such that for any choice of $\mus>0$ and $\muS>0$ one can find
  $\ys=\ys(\mus,\muS)>\frac{1}{R^n}$ and $\gams=\gams(\mus)>0$ with the property that if $T>0$ and $y\in C^1([0,T))$ are
  such that
  \be{5.1}
	y(t) \ge \ys
	\qquad \mbox{for all } t\in (0,T)
  \ee
  and
  \be{5.2}
	y'(t) \le \gams y^{1+\dels}(t)
	\qquad \mbox{for all } t\in (0,T),
  \ee
  then whenever $\theta>0$,
  the functions $\uU=\uU^{(\mus,\al,\theta,y)}$ and
  $\uW=\uW^{(\mus,\beta,\theta,y)}$ from (\ref{uU}) and (\ref{uW}) satisfy
  \be{5.3}
	\parab^{(\muS)} [\uU,\uW](s,t) \le 0
	\qquad \mbox{for all $t\in (0,T)\cap (0,\frac{1}{\theta})$ and $s\in (0,\frac{1}{y(t)})$}.
  \ee
\end{lem}
\proof
  We first fix
  \be{5.4}
	\dels:=1-\beta \in (0,1).
  \ee
  Given $\mus>0$, $\muS>0$ and with $a=a^{(\mus)}$ taken from
  (\ref{a}), we then choose $\ys=\ys(\mus,\muS)>\frac{1}{R^n}$
  large enough such that
  \be{5.5}
	\ys \ge 1
	\qquad \mbox{and} \qquad
	\ys \ge \Big(\frac{2e\muS}{na}\Big)^\frac{1}{1-\beta}
  \ee
  and
  \be{5.6}
	\gams=\gams(\mus):=\frac{nae^{-1}}{2(1-\al)}.
  \ee
  Since  $(e^{-\theta t} \hU)_t = e^{-\theta t} \hU_t - \theta e^{-\theta t} \hU$ for $t\in (0,T)$ and $s\in (0,R^n)$,
  it follows from (\ref{P}) that
  \bea{5.7}
	& & \hs{-20mm}
	e^{\theta t} \cdot \parab^{(\muS)} [\uU,\uW](s,t) \nn\\
	&=& \hU_t - \theta \hU - n^2 s^{2-\frac{2}{n}}  \hU_{ss}
	- ne^{-\theta t} \hU_s \cdot \Big(\hW-\frac{\muS e^{\theta t} s}{n} \Big)
	 \nn\\
	&\le& \hU_t
	- ne^{-\theta t} \hU_s \cdot \Big(\hW-\frac{\muS e^{\theta t} s}{n} \Big)
	\qquad \mbox{for all $t\in (0,T)$ and $s\in (0,\frac{1}{y(t)})$,}
  \eea
  because $\hU \ge 0$ and because  $\hU_{ss} \equiv 0$ in the considered region according to (\ref{1.4}).
  Here we note that due to (\ref{hW}), (\ref{5.5}), (\ref{5.1}) and our assumption $\beta<1$,
  \bea{5.8}
  \hspace{-5mm}
  	\hW - \frac{\muS e^{\theta t} s}{n}
	\ge \hW - \frac{\muS es}{n}
    &= & a y^{1-\beta}(t) s - \frac{\muS es}{n} \nn\\
    &\ge&  a y^{1-\beta}(t) s -\frac{1}{2} a\ys^{1-\beta} s
    \ge  a y^{1-\beta}(t) s -\frac{1}{2} ay^{1-\beta}(t) s
    = \frac{1}{2} a y^{1-\beta}(t) s \nn\\
	& & \qquad \mbox{for all $t\in (0,T)\cap (0,\frac{1}{\theta})$ and $s\in (0,\frac{1}{y(t)})$}
  \eea
  and that by (\ref{1.3}),
  \be{5.9}
  ne^{-\theta t} \hU_s = ne^{-\theta t} \cdot a y^{1-\al}(t)
  \ge nae^{-1} y^{1-\al}(t)
  \qquad \mbox{for all $t\in (0,T)\cap (0,\frac{1}{\theta})$ and $s\in (0,\frac{1}{y(t)})$},
  \ee
  whence using (\ref{1.2}), (\ref{5.8}) and (\ref{5.9})  we infer from (\ref{5.7})
  together with (\ref{5.6}) and (\ref{5.4}) that
  for all $t\in (0,T)\cap (0,\frac{1}{\theta})$ and $s\in (0,\frac{1}{y(t)})$,
    \bas
	e^{\theta t} \cdot \parab^{(\muS)} [\uU,\uW](s,t)
	&\le& (1-\al) a y^{-\al}(t) y'(t) s
	-  nae^{-1} y^{1-\al}(t) \cdot \frac{a}{2} y^{1-\beta}(t) s \\\
	&=& (1-\al)a y^{-\al}(t) s \cdot \bigg\{ y'(t) - \frac{nae^{-1}}{2(1-\al)} y^{2-\beta}(t)\bigg\} \\
    &=& (1-\al)a y^{-\al}(t) s \cdot \bigg\{ y'(t) - \gams y^{1+\dels}(t)\bigg\}\\
    &\le& 0
    \eas
  due to (\ref{5.2}), and this shows (\ref{5.3}).
\qed
We next show that $\qarab [\uU,\uW](s,t) \le 0$ near $s=0$.
\begin{lem}\label{lem6}
  Let $n\in \{3, 4\}$ and $R>0$, let $\mus>0$ and $\muS>0$, assume (\ref{sig}), and suppose that with $T>0$,
  the function $y\in C^1([0,T))$ satisfies
  \be{6.1}
	y(t)>\max \Big\{ 1 \, , \, \frac{1}{R^n}\Big\}
	\quad \mbox{and} \quad
	0\le y'(t) \le  L y^{1+\frac{2}{n}}(t)
	\qquad \mbox{for all } t\in (0,T)
  \ee where $L: =K(e^{-1}a)^{\sig-1} $ with $K$ defined in (\ref{Q}) and $a$ given by (\ref{a}).
  Then whenever
  \be{6.3}
	\theta\ge 1,
  \ee
  for $\uU=\uU^{(\mus,\al,\theta,y)}$ and
  $\uW=\uW^{(\mus,\beta,\theta,y)}$ as in (\ref{uU}) and (\ref{uW}) we have
  \be{6.4}
	\qarab [\uU,\uW](s,t) \le 0
	\qquad \mbox{for all $t\in (0,T) \cap (0,\frac{1}{\theta})$ and $s\in (0,\frac{1}{y(t)})$.}
  \ee
\end{lem}
\proof Since $\sig>\frac{4}{n}$ by (\ref{sig}) and since $n\in \{3, 4\}$,
   we can apply Lemma \ref{lem30} to see that the left-side inequality in (\ref{30.2})
   entails
   \be{6.5}
   \beta+(1-\al)\sig >1+\frac{2}{n}.
   \ee
  For arbitrary $\theta\ge 1$, we recall (\ref{Q}), Lemma \ref{lem2} and Lemma \ref{lem1} to see that
   by the facts that $\hW_{ss}=0$ and that $\hW\ge 0$ in the region $(0,\frac{1}{y(t)})$,
  \bas
	e^{\theta t} \cdot \qarab [\uU,\uW](s,t)
	&=& \hW_t - \theta \hW - n^2 s^{2-\frac{2}{n}} \hW_{ss} + \hW - e^{\theta t} \cdot K s^{1-\sig} e^{-(\theta t)\sig}\hU^\sig \\
	&\le& \hW_t -  K s^{1-\sig} e^{-(\theta t)(\sig-1)}\hU^\sig  \\
	&\le & (1-\beta) a y^{-\beta}(t) y'(t) s
	       -K s^{1-\sig} e^{-(\sig-1)}\cdot a^\sig y^{(1-\al)\sig}(t) s^\sig\\
    &\le &  a y^{-\beta}(t) s
	     \cdot \Big\{y'(t)  - L y^{\beta +(1-\al)\sig}(t) \Big\} \\
    &\le&  a y^{-\beta}(t) s
	     \cdot \Big\{y'(t)  - L y^{1+\frac{2}{n}}(t) \Big\} \\
    &\le& 0
	\qquad \mbox{for all $t\in (0,T) \cap (0,\frac{1}{\theta})$ and $s\in (0,\frac{1}{y(t)})$}
  \eas due to $\sig>1$, $\beta\in (0, 1)$, (\ref{6.1}) and (\ref{6.5}), and this leads to (\ref{6.4}).
 \qed
\subsection{Subsolution properties: intermediate annulus}\label{sect_annulus}		
For our choice of $\beta$, we also have that $\parab^{(\muS)} [\uU,\uW](s,t) \le 0$
in an intermediate region.
\begin{lem}\label{lem7}
  Let $n\in \{3, 4\}$ and $R>0$, assume (\ref{sig}), take any $\delss\in (0, 1-\beta)$,
  and let $\mus>0$ and $\muS>0$ be given.
  Then it is possible to fix
  $\sst \in (0,1]$ with the property that given any $T>0$ and $y\in C^1([0,T))$ fulfilling
  \be{7.1}
	y(t)>\max \Big\{ 1 \, , \, \frac{1}{R^n}\Big\}
	\quad \mbox{and} \quad
	0\le y'(t) \le  y^{1+\delss}(t)
	\qquad \mbox{for all } t\in (0,T),
  \ee
 for any $\theta>0$,
  with $\uU=\uU^{(\mus,\al,\theta,y)}$ and $\uW=\uW^{(\mus,\beta,\theta,y)}$ taken from (\ref{uU}) and (\ref{uW}) we have
  \be{7.2}
	\parab^{(\muS)} [\uU,\uW](s,t) \le 0
	\qquad \mbox{for all $t\in (0,T) \cap (0,\frac{1}{\theta})$ and $s\in (\frac{1}{y(t)},R^n) \cap (0,\sst)$.}
  \ee
\end{lem}
\proof By (\ref{sig}) and Lemma \ref{lem30} we first see that
 \be{7.3}
 1-\frac{2}{n}-\beta>0
 \ee due to $n\in \{3, 4\}$. We then can choose $\sst \in (0,1]$
 sufficiently small such that
  \be{7.4}
  \sst\le \Big(\frac{na}{2e\muS}\Big)^\frac{1}{1-\beta}
  \ee
  and
  \be{7.5}
	\frac{2n(1-\al)e}{a\al} \cdot \sst^{1-\frac{2}{n}-\beta}\le \frac{1}{2}
  \ee
  as well as
  \be{7.6}
	\frac{2(1-\al)e}{na}
		\cdot \sst^{1-\beta-\delss}
	\le \frac{1}{2}.
  \ee due to (\ref{7.3}) and the fact that $1-\beta-\delss>0$ by our choice of $\delss$.
  Now we recall the identity
  \bea{7.7}
  \hspace*{-5mm}
	e^{\theta t} \parab^{(\muS)} [\uU,\uW](s,t)
	&=& \hU_t - \theta \hU - n^2 s^{2-\frac{2}{n}} \hU_{ss}
     -n e^{-\theta t} \hU_s \cdot \Big( \hW - \frac{\muS e^{\theta t} s}{n}\Big), \nn\\[1.5mm]
	& & \hspace{47mm} t\in (0,T), \mbox{$s\in (\frac{1}{y(t)},R^n)$.}
  \eea
  Here we note that
  \be{7.8}
	 s-\frac{1-\beta}{y(t)}
	\ge s-(1-\beta) s = \beta s
	\qquad \mbox{for all $t\in (0,T)$ and $s>\frac{1}{y(t)}$,}
  \ee and thus according to this and (\ref{hW}) together with (\ref{7.4}) we see that
  for all $t\in (0,T)\cap (0,\frac{1}{\theta})$ and $s\in (\frac{1}{y(t)},R^n) \cap (0,\sst)$.
  \bea{7.9}
	\hW - \frac{\muS e^{\theta t} s}{n}
	\ge \beta^{-\beta}a \cdot (\beta s)^\beta - \frac{\muS e^{\theta t} s}{n}
	&\ge& a s^\beta - \frac{\muS e }{n}\cdot s \nn\\
	& =& \frac{1}{2} a s^\beta + \frac{1}{2} a s^\beta \cdot \Big\{1 -\frac{2e\muS  }{na} \cdot s^{1-\beta}\Big\}\nn\\
    &\ge& \frac{1}{2} a s^\beta + \frac{1}{2} a s^\beta \cdot \Big\{1 -\frac{2e\muS  }{na} \cdot \sst^{1-\beta}\Big\}\nn\\
    &\ge& \frac{1}{2} a s^\beta + \frac{1}{2} a s^\beta \cdot \Big\{1 -\frac{2e\muS  }{na} \cdot \frac{na}{2e\muS}\Big\}\nn\\
    &=& \frac{1}{2} a s^\beta.
   %\qquad \mbox{for all $t\in (0,T)\cap (0,\frac{1}{\theta})$ and $s\in (\frac{1}{y(t)},R^n) \cap (0,\ssst)$}.
  \eea Thus, due to the nonnegativity of $\hU$ and $\hU_s$, the identity in
  (\ref{7.7}) can be turned into the inequality
  \bea{7.10}
	e^{\theta t} \parab^{(\muS)} [\uU,\uW](s,t)
	&\le& \hU_t  - n^2 s^{2-\frac{2}{n}} \hU_{ss}
	 -n e^{-\theta t} \hU_s \cdot \frac{a}{2} s^\beta\nn\\[2mm]
	& & \hs{10mm}
	\mbox{for all $t\in (0,T)\cap (0,\frac{1}{\theta})$ and $s\in (\frac{1}{y(t)},R^n) \cap (0,\sst)$}.
  \eea
  Here by (\ref{1.3}) and (\ref{1.4}) and in view of the inequality
  \be{7.11}
  s-\frac{1-\al}{y(t)}
	\ge s-(1-\al) s = \al s
\qquad \mbox{for all $t\in (0,T)$ and $s>\frac{1}{y(t)}$},
\ee
 we have
  \bea{7.12}
	\frac{- n^2 s^{2-\frac{2}{n}} \hU_{ss}}
		{n e^{-\theta t} \hU_s \cdot \frac{a}{2} s^\beta}
	&=& - \frac{2ne^{\theta t}}{a} \cdot s^{2-\frac{2}{n}-\beta} \hU_s^{-1} \hU_{ss} \nn\\
	&=& \frac{2ne^{\theta t}}{a} \cdot s^{2-\frac{2}{n}-\beta} \cdot
		\Big\{ \al^{1-\al} a \cdot \Big(s-\frac{1-\al}{y(t)}\Big)^{\al-1} \Big\}^{-1} \times
		\al^{1-\al} (1-\al) a \cdot \Big(s-\frac{1-\al}{y(t)}\Big)^{\al-2} \nn\\
	&=& \frac{2n(1-\al)e^{\theta t}}{a} \cdot s^{2-\frac{2}{n}-\beta}
			\cdot \Big(s-\frac{1-\al}{y(t)}\Big)^{-1} \nn\\
    &\le& \frac{2n(1-\al)e}{a \al} \cdot s^{1-\frac{2}{n}-\beta} \nn\\
    &\le& \frac{2n(1-\al)e}{a \al} \cdot \sst^{1-\frac{2}{n}-\beta} \nn\\
	&\le & \frac{1}{2}
	\qquad \mbox{for all $t\in (0,T) \cap (0,\frac{1}{\theta})$ and $s\in (\frac{1}{y(t)},R^n) \cap (0,\sst)$}
  \eea due to (\ref{7.5}), which implies the inequality
  \bea{7.13}
	- n^2 s^{2-\frac{2}{n}} \hU_{ss}
	\le \frac{1}{2} \cdot n e^{-\theta t} \hU_s \cdot \frac{a}{2} s^\beta
	\qquad	\mbox{for all $t\in (0,T)\cap (0,\frac{1}{\theta})$ and $s\in (\frac{1}{y(t)},R^n) \cap (0,\sst)$}.
  \eea
  We next use (\ref{1.2}) and (\ref{1.3}) along with (\ref{7.11}) and (\ref{7.1}) to compute
  \bea{7.14}
	\frac{\hU_t}{n e^{-\theta t} \hU_s \cdot \frac{a}{2} s^\beta}
	&=& \frac{2e^{\theta t}}{na}\cdot s^{-\beta} \hU_s^{-1} \hU_t \nn\\
	&=& \frac{2e^{\theta t}}{na} \cdot s^{-\beta} \cdot
		\Big\{ \al^{1-\al} a\cdot \Big(s-\frac{1-\al}{y(t)} \Big)^{\al-1} \Big\}^{-1} \times
    % \nn\\
	% & & \hs{50mm}
		\al^{1-\al} (1-\al) a\cdot \Big(s-\frac{1-\al}{y(t)}\Big)^{\al-1} \cdot \frac{y'(t)}{y^2(t)} \nn\\
	&=& \frac{2(1-\al)e^{\theta t}}{na} \cdot s^{-\beta} \cdot  \frac{y'(t)}{y^2(t)} \nn\\
   &\le& \frac{2(1-\al)e}{na} \cdot s^{-\beta} \cdot \Big(\frac{1}{y(t)}\Big)^{1-\delss} \nn \\
   &\le& \frac{2(1-\al)e}{na} \cdot s^{1-\beta-\delss} \nn\\
   &\le& \frac{2(1-\al)e}{na} \cdot \sst^{1-\beta-\delss} \nn\\
	&\le& \frac{1}{2}
	\qquad \mbox{for all $t\in (0,T)$ and $s\in (\frac{1}{y(t)},R^n) \cap (0,\sst)$}
  \eea
  thanks to (\ref{7.6}), which leads to
    \bas
	\hU_t
	\le \frac{1}{2} \cdot n e^{-\theta t} \hU_s \cdot \frac{a}{2} s^\beta
	\qquad \mbox{for all $t\in (0,T)\cap (0,\frac{1}{\theta})$ and $s\in (\frac{1}{y(t)},R^n) \cap (0,\sst)$}.
  \eas
 Inserting this and (\ref{7.13})  into (\ref{7.10}) establishes (\ref{7.2}).
\qed
Similarly, relying on our choice of the parameters $\alpha$ and $\beta$, we can
find suitable parameter function $y(t)$ and parameter $\theta$
such that $\qarab [\uU,\uW](s,t) \le 0$ in an appropriate intermediate annulus.
\begin{lem}\label{lem9}
  Let $n\in \{3, 4\}$ and $R>0$, assume (\ref{sig}),
  and let $\mus>0$ and $\muS>0$ be given. Then one can find
  $\ssst\in (0,1]$ with the property that if $T>0$ and $y\in C^1([0,T))$ are such that
  \be{9.1}
	y(t)>\max \Big\{ 1 \, , \, \frac{1}{R^n}\Big\}
	\quad \mbox{and} \quad
	0\le y'(t) \le y^{1+\frac{2}{n}}(t)
	\qquad \mbox{for all } t\in (0,T),
  \ee
  and if
  \be{9.2}
	\theta\ge 1,
  \ee
  then the functions $\uU=\uU^{(\mus,\al,\theta,y)}$ and
  $\uW=\uW^{(\mus,\beta,\theta,y)}$ defined in (\ref{uU}) and (\ref{uW}) satisfy
  \be{9.3}
	\qarab [\uU,\uW](s,t) \le 0
	\qquad \mbox{for all $t\in (0,T) \cap (0,\frac{1}{\theta})$ and $s\in (\frac{1}{y(t)},R^n) \cap (0,\ssst)$.}
  \ee
\end{lem}
\proof
  Since
  $\beta-\frac{2}{n}-1 +\sig -\al \sig>0$ asserted by Lemma \ref{lem30}, we can choose
  $\ssst \in (0,1]$ such that
  \be{9.4}
	a(1+n^2 \beta^{-1}) \cdot \ssst^{\beta-\frac{2}{n}-1 +\sig -\al \sig} \le K a^\sig e^{-(\sig-1)}
  \ee with $a$ as in (\ref{a}).
  Due to the nonnegativity of $\hW$ and thanks to $\sig>1$, we can utilize (\ref{Q}), (\ref{9.2}), (\ref{2.2}), (\ref{2.4}) and (\ref{hU})
  together with (\ref{7.8}), (\ref{7.11}) and (\ref{9.1}) to estimate
   \bea{9.5}
   \hspace*{-8mm}
   e^{\theta t} \cdot \qarab [\uU,\uW](s,t)
	&=& \hW_t - \theta \hW - n^2 s^{2-\frac{2}{n}} \hW_{ss} + \hW - e^{\theta t} \cdot K s^{1-\sig} e^{-(\theta t)\sig}\hU^\sig  \nn\\
    &\le& \hW_t - n^2 s^{2-\frac{2}{n}} \hW_{ss}  -K e^{-\theta t(\sig-1)} \cdot s^{1-\sig} \hU^\sig \nn \\
	&=& \beta^{1-\beta} (1-\beta) a \cdot \Big( s-\frac{1-\beta}{y(t)}\Big)^{\beta-1} \cdot \frac{y'(t)}{y^2(t)} \nn\\
	& & + n^2 s^{2-\frac{2}{n}} \cdot \beta^{1-\beta} (1-\beta) a \cdot \Big( s-\frac{1-\beta}{y(t)}\Big)^{\beta-2} \nn\\
	& & -K e^{-\theta t(\sig-1)} \cdot s^{1-\sig} \cdot \Big\{ \al^{-\al} a\cdot \Big( s- \frac{1-\al}{y(t)} \Big)^\al\Big\}^\sig \nn\\
	&\le& \beta^{1-\beta} (1-\beta) a \cdot (\beta s)^{\beta-1} \cdot y^{-1+\frac{2}{n}}(t) \nn\\
	& & + n^2 s^{2-\frac{2}{n}} \cdot \beta^{1-\beta} (1-\beta)a \cdot  (\beta s)^{\beta-2} \nn\\
	& & - K e^{-(\sig-1)} \cdot s^{1-\sig} \cdot \{ \al^{-\al} a \cdot (\al s)^\al\}^\sig \nn\\
    &\le& a(1+n^2 \beta^{-1}) \cdot s^{\beta -\frac{2}{n}} - K a^\sig e^{-(\sig-1)} \cdot s^{1-\sig+\al\sig} \nn\\
    &=&  s^{1-\sig+\al\sig} \cdot \Big\{ a(1+n^2 \beta^{-1}) \cdot s^{\beta -\frac{2}{n}-1 +\sig -\al\sig}
       -K a^\sig e^{-(\sig-1)} \Big\} \nn\\
    &=&  s^{1-\sig+\al\sig} \cdot \Big\{ a(1+n^2 \beta^{-1}) \cdot \ssst^{\beta -\frac{2}{n}-1 +\sig -\al\sig}
       -K a^\sig e^{-(\sig-1)} \Big\} \nn\\
	&\le& 0
   \qquad \mbox{for all $t\in (0,T) \cap (0,\frac{1}{\theta})$ and $s\in (\frac{1}{y(t)},R^n) \cap (0,\ssst)$},
  \eea due to our choice of $\ssst$ in (\ref{9.4}), and this yields (\ref{9.3}).
 \qed

\subsection{Subsolution properties: outer region}\label{sect_outer}		
Depending on the properties of the function $\hU$ in an outer region, we can find some large
parameter $\theta$ and parameter function $y(t)$ such that $\parab^{(\muS)} [\uU,\uW](s,t) \le 0$
holds in the outer region.
\begin{lem}\label{lem8}
  Let $n\in \{3, 4\}$, $R>0$, $\mus>0$ and $\muS>0$, assume (\ref{sig}),
  and
  let $s_0\in (0,1]$ be arbitrarily given.
  Then there exists $\ths=\ths^{(s_0)}>0$ such that whenever $\theta \ge \ths$ and $T>0$ as well as $y\in C^1([0,T))$
  are such that
  \be{8.1}
	y(t)>\frac{1}{R^n}
	\quad \mbox{and} \quad
	 0\le y'(t) \le y^2(t)
	\qquad \mbox{for all } t\in (0,T),
  \ee
 the functions $\uU$ and $\uW$ defined in (\ref{uU}) and (\ref{uW}) satisfy
  \bea{8.2}
	& & \hs{-32mm}
	\parab^{(\muS)} [\uU,\uW](s,t) \le 0
	\qquad \mbox{for all $t\in (0,T) \cap (0,\frac{1}{\theta})$ and $s\in (\frac{1}{y(t)},R^n) \cap [s_0,R^n)$.}
  \eea
\end{lem}
\proof
  For any given $s_0\in (0, 1]$, we first choose $\ths=\ths^{(s_0)}>0$ sufficiently large such that
  \be{8.3}
	\ths > \frac{1+\muS R^n}{s_0^{1+\al}} +\frac{n^2 R^{2n-2}}{\al s_0^{2+\al}}.
  \ee
  Since
  \bas
	s-\frac{1-\al}{y(t)} \ge s-(1-\al)s=\al s
	\ge \al s_0
	\qquad \mbox{for all $t\in (0,T)$ and $s\in (\frac{1}{y(t)},R^n) \cap [s_0,R^n)$}
  \eas
  and since $\al\in (0, 1)$ and $s_0\le 1$,
  we then obtain that according to (\ref{hU}), (\ref{1.2}), (\ref{1.3}) and (\ref{1.4})
  along with (\ref{8.1}),
  \bas
	\hU = \al^{-\al} a \cdot \Big(s-\frac{1-\al}{y(t)}\Big)^\al
	\ge \al^{-\al} a\cdot (\al s_0)^\al
	= a s_0^\al
  \eas
  and
  \bas
	\hU_t
	&=& \al^{1-\al} (1-\al) a \cdot \Big(s-\frac{1-\al}{y(t)}\Big)^{\al-1} \cdot \frac{y'(t)}{y^2(t)} \\
	&\le& \al^{1-\al} (1-\al) a \cdot (\al s_0)^{\al-1}
	= (1-\al) a s_0^{\al-1}
	\le \frac{a}{s_0}
  \eas
  as well as
  \bas
	\hU_s
	= \al^{1-\al} a \cdot \Big(s-\frac{1-\al}{y(t)}\Big)^{\al-1}
	\le \al^{1-\al} a \cdot (\al s_0)^{\al-1}
	= a s_0^{\al-1}
	\le \frac{a}{s_0}
  \eas
  and
  \bas
	-\hU_{ss} = \al^{1-\al} (1-\al) a \cdot \Big(s-\frac{1-\al}{y(t)}\Big)^{\al-2}
	\le \al^{1-\al} (1-\al) a \cdot (\al s_0)^{\al-2}
	= \frac{(1-\al) a s_0^{\al-2}}{\al}
	\le \frac{a}{\al s_0^2}
  \eas
  for all $t\in (0,T)$ and $s\in (\frac{1}{y(t)},R^n) \cap [s_0,R^n)$.
  Relying on the above estimates and utilizing the nonnegativity of $\hU_s$ and $\hW$,
  we see that whenever $\theta>\ths$,
  \bas
	\hspace*{-5mm}
	e^{\theta t} \parab^{(\muS)} [\uU,\uW](s,t)
	&=& \hU_t - \theta \hU - n^2 s^{2-\frac{2}{n}} \hU_{ss}
       -n e^{-\theta t} \hU_s \cdot \Big( \hW - \frac{\muS e^{\theta t} s}{n}\Big) \nn\\
    &\le & \hU_t - \theta \hU - n^2 s^{2-\frac{2}{n}} \hU_{ss} +\muS s \cdot \hU_s  \nn\\
    &\le& \frac{a}{s_0} -\theta a s_0^\al +n^2 (R^n)^{2-\frac{2}{n}} \cdot \frac{a}{\al s_0^2} +\muS R^n \cdot \frac{a}{s_0}\nn\\
    &=& a s_0^\al \cdot \Big\{\frac{1+\muS R^n}{s_0^{1+\al}} +\frac{n^2 R^{2n-2}}{\al s_0^{2+\al}} -\theta\Big\} \nn\\
    &\le& a s_0^\al \cdot \Big\{\frac{1+\muS R^n}{s_0^{1+\al}} +\frac{n^2 R^{2n-2}}{\al s_0^{2+\al}} -\ths\Big\} \nn\\
	&\le & 0
    \qquad \mbox{for all $t\in (0,T) \cap (0,\frac{1}{\theta})$ and $s\in (\frac{1}{y(t)},R^n) \cap [s_0,R^n)$}
  \eas thanks to (\ref{8.3}). This yields (\ref{8.2}).
 \qed
Similarly, using several explicit estimates of the function $\hW$ in an outer region, we can chose sufficiently large
parameter $\theta$ and parameter function $y(t)$ such that $\qarab [\uU,\uW](s,t) \le 0$ is satisfied
in the outer region.
\begin{lem}\label{lem10}
  Let $n\in \{3, 4\}$, $R>0$, $\mus>0$ and $\muS>0$,  and assume (\ref{sig}).
   Then for any given $s_0\in (0,1]$ there exists $\thss=\thss^{(s_0)}>0$ such that if
 $\theta \ge \thss$,  $T>0$ and $y\in C^1([0,T))$ are such that
  \be{10.1}
	y(t)>\max\Big\{ 1\, , \, \frac{1}{R^n}\Big\}
	\quad \mbox{and} \quad
	 0\le y'(t) \le y^{1 +\frac{2}{n}}(t)
	\qquad \mbox{for all } t\in (0,T),
  \ee
  with $\uU=\uU^{(\mus,\al,\theta,y)}$ and
  $\uW=\uW^{(\mus,\beta,\theta,y)}$ taken from (\ref{uU}) and (\ref{uW}), there holds
  \be{10.2}
	\qarab [\uU,\uW](s,t) \le 0
	\qquad \mbox{for all $t\in (0,T)$ and $s\in (\frac{1}{y(t)},R^n) \cap [s_0,R^n)$.}
  \ee
\end{lem}
\proof
  Given any $s_0\in (0, 1]$, we take $\thss=\thss^{(s_0)}>0$ large enough such that
  \be{10.3}
	\thss\ge 2
	\qquad \mbox{and} \qquad
	\thss \ge2(1-\beta)(1+n^2 \beta^{-1}) s_0^{-\frac{n}{2}}.
  \ee
  Again since
  \bas
  s-\frac{1-\beta}{y(t)} \ge \beta s
  \quad\mbox{and}\quad
  s-\frac{1-\al}{y(t)}\ge \al s
  \qquad\mbox{$t\in (0,T)$ and $s\in (\frac{1}{y(t)},R^n)$},
   \eas
 in view of  Lemma \ref{lem2} together with the nonnegativity of $\hU$ and (\ref{10.1}),
  we obtain that if $\theta \ge \thss$,
   \bea{10.4}
   \hspace{-5mm}
	 e^{\theta t} \cdot \qarab [\uU,\uW](s,t)
	&=& \hW_t - \theta \hW - n^2 s^{2-\frac{2}{n}} \hW_{ss} + \hW - e^{\theta t} \cdot K s^{1-\sig} e^{-(\theta t)\sig}\hU^\sig  \nn\\
    &\le& \hW_t - n^2 s^{2-\frac{2}{n}} \hW_{ss} -\frac{\theta}{2} \hW  \nn\\
    &=& \beta^{1-\beta} (1-\beta) a \cdot \Big( s-\frac{1-\beta}{y(t)}\Big)^{\beta-1} \cdot \frac{y'(t)}{y^2(t)} \nn\\
	& & + n^2 s^{2-\frac{2}{n}} \cdot \beta^{1-\beta} (1-\beta) a \cdot \Big( s-\frac{1-\beta}{y(t)}\Big)^{\beta-2} \nn\\
    & & -\frac{\theta}{2} \beta^{-\beta} a \cdot  \Big( s-\frac{1-\beta}{y(t)}\Big)^{\beta} \nn\\
    &\le& \beta^{1-\beta} (1-\beta) a \cdot (\beta s)^{\beta-1} \cdot \Big(\frac{1}{y(t)}\Big)^{1-\frac{2}{n}} \nn\\
    & & + n^2 s^{2-\frac{2}{n}} \cdot \beta^{1-\beta} (1-\beta) a \cdot (\beta s)^{\beta -2} \nn\\
    & &  -\frac{\theta}{2} \beta^{-\beta} a \cdot (\beta s)^\beta \nn\\
    &\le& \beta^{1-\beta} (1-\beta) a \cdot (\beta s)^{\beta-1} \cdot s^{1-\frac{2}{n}}
           +n^2\beta^{-1} (1-\beta) a \cdot s^{\beta -\frac{2}{n}} -\theta \cdot \frac{a}{2} s^\beta \nn\\
    &=& (1-\beta)(1+n^2 \beta^{-1}) a \cdot s^{\beta -\frac{2}{n}} -\theta \cdot \frac{a}{2} s^\beta \nn\\
    &=&  \frac{a}{2} s^\beta  \cdot \Big\{2(1-\beta)(1+n^2 \beta^{-1}) a \cdot s^{-\frac{2}{n}} -\theta\Big\} \nn\\
    &\le&  \frac{a}{2} s^\beta  \cdot \Big\{2(1-\beta)(1+n^2 \beta^{-1}) a \cdot s_0^{-\frac{2}{n}} -\thss\Big\} \nn\\
	&\le& 0	\qquad \mbox{for all $t\in (0,T)$ and $s\in (\frac{1}{y(t)},R^n) \cap [s_0,R^n)$,}
  \eea
  because $0<\beta<1$ and because (\ref{10.3}). The completes the proof.
 \qed
\subsection{Proof of Theorem \ref{theo13}}\label{sect_BU_proof}
We now make sure that under the assumption on $\sig$ from Theorem \ref{theo13} we can chose the parameters
$\alpha, \beta, \theta$ and the parameter function $y(t)$ such that all hypotheses in Lemmata \ref{lem5}-\ref{lem10}
are simultaneously satisfied, and thus the pair $(\uU,\uW)$ will indeed forms a subsolution.
\begin{lem}\label{lem12}
   Let $n\in \{3, 4\}$, $R>0$, $\mus>0$ and $\muS>0$, and assume (\ref{sig}).
   Then for any $\Ts>0$,
  one can find $\theta>0, T\in (0,\Ts)$ and a positive function $y\in C^1([0,T))$ such that
  \be{12.1}
	y(t) \to +\infty
	\qquad \mbox{as } t\nearrow T,
  \ee
  and that with the functions $\uU=\uU^{(\mus,\al,\theta,y)}$ and
  $\uW=\uW^{(\mus,\beta,\theta,y)}$ taken from (\ref{uU}) and (\ref{uW}) satisfy
  \be{12.2}
	\parab^{(\muS)} [\uU,\uW](s,t) \le 0
	\quad \mbox{and} \quad
	\qarab [\uU,\uW](s,t) \le 0
	\qquad \mbox{for all $t\in (0,T)$ and $s\in (0,R^n)\sm \{ \frac{1}{y(t)} \}$.}
  \ee
\end{lem}
\proof
  With $\dels>0$ and $\delss>0$ provided by Lemma \ref{lem5} and Lemma \ref{lem7}, we take
  \be{12.3}
	\del:=\min \Big\{ \dels \, , \, \delss \, , \, \frac{2}{n} \Big\}.
  \ee
  We then let $\sst=\sst^{(\delss)}$ and $\ssst$
  be chosen in Lemma \ref{lem7} and Lemma \ref{lem9}, and apply Lemma \ref{lem8} and Lemma \ref{lem10} to
  \be{12.4}
	s_0:=\min \big\{ \sst \, , \, \ssst \big\}
  \ee
  to obtain $\ths=\ths^{(s_0)}$ and $\thss=\thss^{(s_0)}$ with the properties therein. Setting
  \be{12.5}
	\theta:=\max \big\{ \ths \, , \, \thss \, , \, 2 \big\}
  \ee
  and taking $\gams$, $\ys$ and $L$
  from lemma \ref{lem5} and Lemma \ref{lem6}, we then let
  \be{12.6}
	\gamma := \min \Big\{ \gams \, , \, L \, ,\, 1\Big\}
  \ee
  and select $y_0>0$ large enough such that
  \be{12.7}
	y_0>\max \Big\{ 1 \, \frac{1}{R^n}\Big\}
	\qquad \mbox{and} \qquad
	y_0 \ge \ys,
  \ee
  and that
  \be{12.8}
	T:=\frac{1}{\gamma\del y_0^\del}
	\quad \mbox{satisfies} \quad
	T < \min \Big\{ \frac{1}{\theta} \, , \, \Ts \Big\}.
  \ee
  Our choice of $T$ guarantees that the problem
  \be{12.9}
	\lbal
	y'(t) = \gamma y^{1+\del}(t),
	\qquad t\in (0,T), \\[1mm]
	y(0)=y_0,
	\ear
  \ee
  has a solution $y\in C^1([0,T))$ which satisfies (\ref{12.1}) as well as $y\ge y_0$, so that (\ref{12.7}),
  (\ref{12.6}) and (\ref{12.3}) ensure that all the requirements on $y$ made in Lemmata \ref{lem5}-\ref{lem10} are fulfilled.
  Moreover, (\ref{12.5}) makes sure that also the corresponding assumptions on $\theta$ in these lemmata
  are met, we may conclude from the above-said statements and from (\ref{12.5}) that for
  $(\uU,\uW)=(\uU^{(\mus,\al,\theta,y)},\uW^{(\mus,\beta,\theta,y)})$
  defined in  (\ref{uU}) and (\ref{uW}) we have
  \bas
	\parab^{(\muS)} [\uU,\uW](s,t)\le 0
	\quad \mbox{and} \quad
	\qarab [\uU,\uW](s,t)\le 0
	\qquad \mbox{for all $t\in (0,T) \cap (0,\frac{1}{\theta})$ and $s\in (0,R^n)\sm \{ \frac{1}{y(t)}\}$.}
  \eas
  Since from (\ref{12.8}) we also have that $(0,T)\cap (0,\frac{1}{\theta})=(0,T)$ and that $T<\Ts$, this completes the proof.
\qed
We are now in a position to complete the proof of Theorem \ref{theo13}.\abs
\proofc of Theorem \ref{theo13}. \quad
  On the basis of Lemma \ref{lem12}, one readily verifies that $(\uU, \, \uW)$ defined through (\ref{uU}) and (\ref{uW})
  indeed forms a sub-solution of the  parabolic system (\ref{11.3})-(\ref{11.4}) and (\ref{11.2}), and thus
  in view of (\ref{uU}) and (\ref{1.3}) we have that
  $u(0, t) =n U_s(0, t)\ge n\uU_s(0, t)=e^{-\theta t}\cdot ay^{1-\al}(t)$, where the latter blows up at $T$ due to (\ref{12.1})
   and $\al<1$.
  The readers may refer \cite{taowin_291} for a detailed proof.
\qed
{\bf Acknowledgement.} \quad
 This work was supported by the {\em National Natural Science Foundation of China
   (No. 12171316)}.
\end{document}